\documentclass[12pt]{article}
\usepackage{color}
\usepackage{graphicx}
\usepackage[colorlinks=true]{hyperref}
\usepackage{amssymb}
\usepackage{amsmath}
\usepackage{amsthm}
\usepackage{mathtools}
\usepackage{authblk}
\usepackage{fullpage}

\theoremstyle{plain}
\newtheorem{theorem}{Theorem}[section]
\newtheorem{lemma}[theorem]{Lemma}

\newtheorem{corollary}[theorem]{Corollary}

\theoremstyle{definition}
\newtheorem{definition}[theorem]{Definition}

\newtheorem{remark}[theorem]{Remark}
\newtheorem{example}[theorem]{Example}

\makeatletter
\@namedef{subjclassname@2010}{%
  \textup{2010} Mathematics Subject Classification}
  \makeatother

\newcommand{\diffint}[5]{{}^{#2}_{#5}\!#1^{(#4)}_{#3}}
\newcommand{\RLI}[2]{\diffint{I}{RL}{#1}{#2}{}}     % riemann Liouville integral

\newcommand{\RL}[2]{\diffint{D}{RL}{#1}{#2}{}}    % riemann Liouville derivative
\newcommand{\CA}[2]{\diffint{D}{C}{#1}{#2}{}}     % Caputo
\newcommand{\RLs}[2]{\diffint{D}{RL}{#1\,*}{#2}{}}  
\newcommand{\CAs}[2]{\diffint{D}{C}{#1\,*}{#2}{}}

\newcommand{\purge}[1]{}

\newcommand{\R}{\mathbb{R}} \newcommand{\Z}{\mathbb{Z}}
\newcommand{\C}{\mathbb{C}} \newcommand{\N}{\mathbb{N}}

\newcommand{\norm}[2]{\lVert#1\rVert_{#2}}
\newcommand{\nheight}[1]{\norm{#1}{\infty}}
\newcommand{\nlength}[1]{\norm{#1}{1}}
\newcommand{\ntwo}[1]{\norm{#1}{2}}
\newcommand{\ff}[2]{{#1}_{#2}}
\DeclareMathOperator{\sign}{sgn} % sign

\allowdisplaybreaks

\title{Zeros of fractional derivatives of polynomials}
\author{Torre Caparatta\thanks{\texttt{tecapara@uncg.edu}}}
\author{Sebastian Pauli\thanks{\texttt{s\_pauli@uncg.edu}}}
\author{Filip Saidak\thanks{\texttt{f\_saidak@uncg.edu}}}
\affil{Department of Mathematics and Statistics, University of North Carolina Greensboro}
\date{\today}

\begin{document}

\maketitle

\begin{abstract}
We investigate the behavior of fractional derivatives of polynomials.  In particular,
we consider the locations and the asymptotic behavior of their zeros and give bounds for their Mahler measure.
\end{abstract}

\section{Introduction}

Questions concerning finding exact or approximate values of the zeros of polynomial functions $p(x) = c_n x^n + c_{n-1}x^{n-1} + \cdots + c_1 x + c_0$ are classical, and (for the case of real coefficients $c_0, c_1, \cdots c_n$) properties of the distribution of these zeros have been studied since at least 1637, when Descartes established his fundamental Rule of Signs (in {\em La G{\' e}om{\' e}trie} \cite{descartes1637}). This important result was refined to finite intervals by Budan in 1807 and by Fourier in 1820 (see \cite{bf1} or \cite{bf2}). By then, thanks to the work of Euler, Gauss and Argand (see \cite{fine1997} for more details), the Fundamental Theorem of Algebra had been established, guaranteeing that a polynomial of degree $n$ has exactly $n$ complex zeros (counted with multiplicity). Not much later,
in 1829, Cauchy \cite{cauchy1829} was able to prove that all zeros of a monic polynomial $p(x)$, with complex coefficients, must lie inside the disk $|z| < 1 + \max_{0 \leq k \leq n-1} |c_k|$, 
the bounds that were eventually generalized by Landau \cite{landau1907}, Fej{\' e}r \cite{fejer1908} and others.  

In another direction,
one could ask about the relation between locations of the zeros of a polynomial $p(x)$ and the zeros of its derivative $p^{\prime}(x)$, as Rolle has done in his {\em  Trait{\' e} d’alg{\`e}bre} of 1690 (see \cite{shain1937}); the well-known theorem bearing his name -- that states that between any two zeros of a real polynomial there lies at least one zero of the derivative -- was proved rigorously by Cauchy \cite{cauchy1823} in 1823. In the complex plane, the situation becomes even more interesting. As Gauss noted in 1836, all zeros of $p^{\prime}(x)$ lie in the convex hull of the zeros of $p(x)$.
The first proof of this proposition was published by Lucas \cite{lucas1874} in 1874; it is now known as the Gauss-Lucas Theorem
(also see \cite{marden1966}). At the beginning of the 20th century it was refined by B{\^o}cher \cite{bocher1904}, Jensen \cite{jensen1912} and Walsh \cite{walsh1920}, and in more recent times, several related extensions and generalizations of it have been considered by Dimitrov \cite{dimitrov1998}, Brown \& Xiang \cite{brown1999}, Sendov \cite{sendov2021}, Tao \cite{tao2022} and others.

\begin{figure}[ht]
\includegraphics[width=0.45\textwidth]{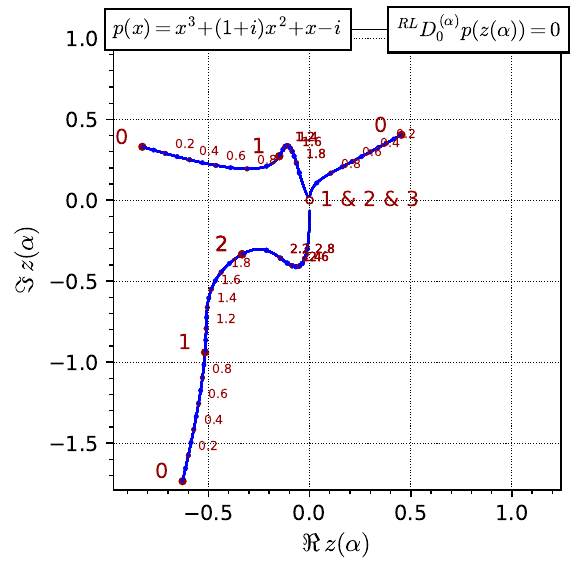}
\includegraphics[width=0.45\textwidth]{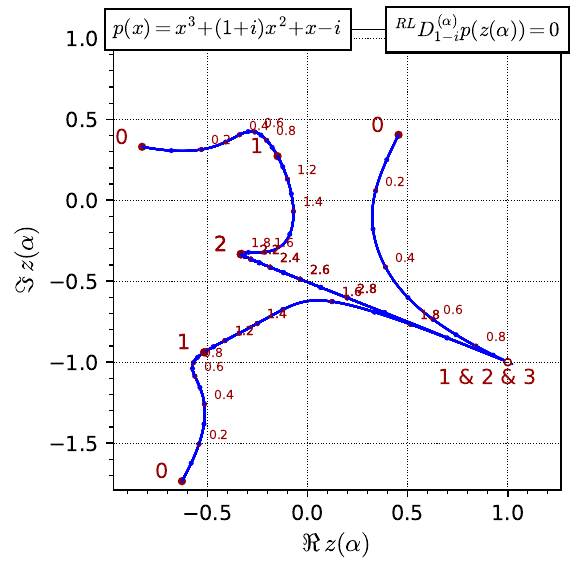}
\caption{Paths \(z(\alpha)\) of zeros of the Riemann-Liouville fractional derivatives 
\(\RL{0}{\alpha}p\) and
\(\RL{1-i}{\alpha}p\)
of the polynomial \(p(x)=x^3 + (1+i)x^2 + x - i\).
We end the paths when \(z(\alpha)\) reaches the ``origin'' (\(a=0\) in the first case, and \(a=1-i\) in the second).
% Also the zeros of integral derivatives that lie on the paths starting at the same zero of \(p(x)\) differ.
}\label{fig two}
\end{figure}

The main aim of our work is to investigate connections between these two central themes. We will try to show that their key ideas can be combined in a very natural way, but to quite surprising effects, if one considers the fractional 
derivatives $p^{(\alpha)} (x)$, where $\alpha \in \mathbb{R}$ is a variable $0 \leq \alpha \leq n = \deg p(x)$.  %Surprisingly, this line of inquiry has not been pursued before. 
Our main goal in this paper will be to answer one of the most intriguing questions that arises as soon as one begins to study these topics: since obviously $\deg p^{\prime}(x) =  \deg p(x) - 1$, and the Fundamental Theorem asserts than the same reduction must occur for the total number of zeros, 
what happens to the zeros of fractional derivatives, as the real $\alpha$ increases continuously from $0$ to $n$? How 
do the zeros of polynomials vanish, and why? As it turns out, these questions have remarkably simple and elegant answers. Namely, for a polynomial $p(x)$ of degree $n$, each of its $n$ zeros will belong to a path of unique length that connects it to the origin, where the ``length'' of the path can be measured by the number of zeros of its derivatives it contains; in other words, for each $0 \leq k \leq n-1$ there will be a unique path (originating at one of the zeros of $p(x)$) that will contain exactly $k$ zeros of its higher derivatives. (Figure \ref{fig two} shows this general property for a generic cubic polynomial).

Another goal of this paper will be to try to understand some of the particulars of the paths the fractional zeros take, their dynamical properties. In order to state our results concerning this general flow of polynomial zeros more precisely, first we will need to recall some basic definitions and properties 
of fractional derivatives. There exist a multitude of different definitions of fractional derivatives, each with its own particular advantages and disadvantages. Unlike in 
our study of the Riemann zeta function and the Stieltjes constants where the Gr\"unwald-Letnikov fractional derivative %$\GL{a}{\alpha} \zeta(s)$
\cite{g:1867,l:1869-1} worked 
really well (see \cite{fps} and \cite{fps2}), in the case of polynomials their divergence forces us instead to employ 
%$\RL{a}{\alpha} p(x)$, 
the  Riemann-Liouville fractional derivative \cite{riefrac}, which also can be thought of as a truncated version of the Gr\"unwald-Letnikov fractional derivative \cite{ortigueira2011fractional}.

The rest of the paper is structured as follows. Our Section \ref{sec derv} contains a short introduction to the Riemann-Liouville differintegral  focusing on results most applicable to fractional derivatives of polynomials. 
In Section \ref{sec low} this theory is applied to the two simplest cases: polynomials of degree one and two. These are the two cases where, thanks to the manageable classical formulas for the zeros, all the main questions can be conclusively answered. With the cubic polynomials things become somewhat murky, but general convergence trends can still be established. In Section \ref{sec flow} we do just that: we examine how the zeros of integral derivatives are connected to the zeros of fractional derivatives in the most general setting, and we look at the paths of zeros and investigate their convergence and the overall flow. In Section \ref{sec bounds} we consider the behaviour of the zeros on a larger scale and we prove bounds for
the Mahler measure of the fractional derivatives, which are then, in Section \ref{caputo}, also established for the Caputo fractional derivative. Finally, in Section \ref{conclusion}, we discuss some intriguing open problems and unsolved questions. 

\section{Riemann-Liouville Fractional 
Derivatives}\label{sec derv}
In full generality, the \(\alpha\)-th Riemann-Liouville fractional derivative is defined as follows, see \cite{oldham1974fractional}, for example: 

\begin{definition}[Riemannn-Liouville fractional derivative]\label{def rl}
Let \(f\) be analytic on a convex open set \(C\) let and \(a\in C\). 
For \(\alpha>0\) the Riemann-Liouville integral is
\[
\RLI{a}{\alpha} f(t)=\frac{1}{\Gamma(\alpha)}
\int_a^t f(\tau)(t-\tau)^{\alpha-1} \,d\tau.
\]
Set \(\RLI{a}{0}f(t)=f(t)\) and for \(\alpha>0\) define the Riemann-Liouville fractional derivative as
\[
\RL{a}{\alpha} f(t) = \frac{d^m}{dt^m} \RLI{a}{m-\alpha} f(t),
\]
where $m=\lceil \alpha \rceil$.  For \(\alpha<0\) set
\(\RL{a}{\alpha}f(t)=\RLI{a}{-\alpha} f(t)\).
\end{definition}
In what follows, we will only consider the special case of polynomials, composed of the simple power functions $p(x)=(x-a)^{\beta}$, where $\beta \in \mathbb{R}$, $a\in\C$. For these, the $\alpha$-th Riemann-Liouville fractional derivative can be computed 
using the Power Rule: 
\begin{align}\label{eq poly derv}
\RL{a}{\alpha}(x-a)^{\beta}=\begin{cases}
0 & \text{ if }\alpha-\beta \in \N \\
%\displaystyle
\dfrac{\Gamma(\beta+1)}{\Gamma(\beta-\alpha+1)}(x-a)^{\beta-\alpha} & \text{otherwise}
\end{cases}
\end{align}

where, for $z \in \mathbb{C}$, the gamma function $\Gamma(z)$ is defined as $\Gamma(z)=\int_0^\infty t^{z-1}e^{-t}dt$; it satisfies $\Gamma(1+n)=n\Gamma(n)$, which implies that $\Gamma(n) = (n-1)!$ for $n \in \mathbb{N}$. The function $\Gamma(z)$ has no zeros in the complex 
plane $\mathbb{C}$, and has poles at all the negative integers (see \cite{gammafunc}).

\begin{remark}
The Riemann-Liouville fractional derivative of a monomial \(f(x)=x^n\) is multivalued.  When changing the branch of the complex logarithm in the computation of the fractional derivative all coefficients of the derivative are changed by the same factor.
So choosing a different branch of the complex logarithm does not change the zeros of the derivative, which means that we
can fix the branch in our consideration of zeros of derivatives of polynomials.
We use the principal branch of the complex logarithm.
\end{remark}

It should be noted that the constant $a$ that centers the expansion (\ref{eq poly derv}) plays a key role in all our computations below, as the ``origin,'' or the limit of convergence, of the flow of zeros of derivatives (Figure \ref{fig two} illustrates its role). Also noteworthy is the fact that the Riemann-Liouville fractional derivative satisfies all properties expected of a regular derivative, with the exception of the composition rule. The following example shows why it fails: 

\begin{example} 
From (\ref{eq poly derv}),  the \(1.5\)-th derivative of \(p(x)=1\) is
\(
\textstyle\RL{0}{1.5}p(x)=\frac{2}{ \sqrt{\pi}}x^{-1.5}
\)
and 
\(\textstyle\RL{0}{1}\left(\RL{0}{0.5}p(x)\right) =\frac{2}{ \sqrt{\pi}}x^{-1.5}.
\)
However, \(\RL{0}{0.5}\left(\RL{0}{1}p(x)\right)=\RL{0}{0.5}0=0.
\)
\end{example}

When \(\beta\in\R\setminus\Z\) we still have
\(\RL{a}{\alpha}\left(\RL{a}{\beta}p(x)\right)=\RL{a}{\alpha+\beta}p(x)\).

\begin{remark} It is possible to go beyond the standard values of $0 \leq \alpha$, and consider what happens for $\alpha < 0$. Here, there will be \(n+1\) complex roots, because the first term in (\ref{eq lem derv}) below has the root \(a\). Just like in the standard case, the extended curves \(z(\alpha)\) of zeros of the differintegral \(\RL{a}{\alpha}p\) will be continuous for \(\alpha<0\) unless
\(\RL{a}{\alpha}p\) has a double root; however, they will not be smooth at integral \(\alpha > 0\).  More on this will be said in Section \ref{sec flow} below.
\end{remark}

In what follows, we consider the zeros of the the fractional derivatives of polynomials \(p\in\C[x]\) of degree \(n\), and we investigate the implicit functions \(z:[0,n)\to \C\) given by
\[
\RL{a}{\alpha}f(z(\alpha))=0.
\]
If \(\left(\RL{a}{\alpha}p(x)\right)'\ne 0\), for \(\alpha\in[0,n)\), then 
\(z(\alpha)\) is differentiable on \([0,n)\).
We denote the roots of the polynomial \(p(x)\) by \(z_1, z_2, \dots,z_n\)
and for \(1\le k\le n\) we define the the implicit function \(z_k:[0,n)\to\C\)  by \(z_k(0)=z_k\)
and \(\RL{a}{\alpha}p(z(\alpha))=0\).

The following representation of the fractional derivatives of a general monic polynomial will be most useful.

\begin{lemma}\label{lem derv}
Let \(p(x)=(x-a)^n+\sum_{j=0}^{n-1} c_j (x-a)^{j}\in\mathbb{C}[x]\) and \(\alpha\in\R\setminus\N^{\ge n}\).  
Set
\begin{equation}\label{eq lem derv}
\RL{a*}{\alpha}p(x):=
(x-a)^{n}+
\sum_{j=j_0}^{n-1} 
\left(
\prod_{k=j+1}^n\!(k-\alpha)\right)\cdot
\frac{j!}{n!} c_j (x-a)^{j}.
\end{equation}
where $j_0=\left\{\begin{array}{ll}\alpha & \text{if } \alpha\in\N\\0 & \text{otherwise}\end{array}\right.$.  Then
\begin{enumerate}
\item If $\RL{a}{\alpha}p(z)=0$ then $\RL{a*}{\alpha} p(z)=0$ or $\alpha<0$ and $z=a$.
\item If $\RL{a*}{\alpha} p(z)=0$ then $\RL{a}{\alpha} p(z)=0$ or $\alpha\in\N$ and $z=a$.
\end{enumerate}
\end{lemma}

\begin{proof}
For $\alpha\not\in\N$
with the Power Rule (\ref{eq poly derv}), we obtain \begin{align*}
%f^{(\alpha)}(x)= 
\RL{a}{\alpha}p(x)  &=
\frac{\Gamma(n+1)}{\Gamma(n+1-\alpha)} (x-a)^{n-\alpha}+
\sum_{j=0}^{n-1} \frac{\Gamma(j+1)}{\Gamma(j+1-\alpha)} c_j (x-a)^{j-\alpha}\\
&=\frac{\Gamma(n+1)}{\Gamma(n+1-\alpha)} (x-a)^{-\alpha}
\left[(x-a)^{n}+
\sum_{j=0}^{n-1} 
\frac{\Gamma(n+1-\alpha)}{\Gamma(j+1-\alpha)}
\frac{\Gamma(j+1)}{\Gamma(n+1)} c_j (x-a)^{j}
\right]
\\
&=\frac{n!}{\Gamma(n+1-\alpha)} (x-a)^{-\alpha}
\left[(x-a)^{n}+
\sum_{j=0}^{n-1} 
\left(
\prod_{k=j+1}^n\!(k-\alpha)\right)\cdot
\frac{j!}{n!} c_j (x-a)^{j}
\right]\\
&=\frac{n!}{\Gamma(n+1-\alpha)} (x-a)^{-\alpha}\cdot\RL{a*}{\alpha}p(x)
\end{align*}
For  $\alpha\in\N^{\le n}$ we have
$\RL{a*}{\alpha)p(x)}=\frac{n!}{(n-\alpha)!}(x-a)^\alpha p^{(\alpha)}(x)$.
\end{proof}

\begin{remark} The representation of the fractional derivatives of a polynomial $p(x)$ given in Lemma \ref{lem derv} has the property that their roots only depend on the factors in $\RL{a*}{\alpha}$, which in turn implies the useful fact that the branch cut of the complex logarithm does not affect the paths \(z(\alpha)\) of
zeros of these fractional derivatives.

We also see that for \(p(x)=(x-a)^n\) we have:
\begin{equation}\label{eq at a}
\RL{a}{\alpha}p(a)=
\begin{cases}
0 &\mbox{if }\alpha<n\\
\Gamma(n+1) & \mbox{if }\alpha=n\\
\mbox{undefined} &\mbox{if }\alpha>n
\end{cases}
\end{equation}
\end{remark}

\begin{figure}
\includegraphics[width=0.45\textwidth]{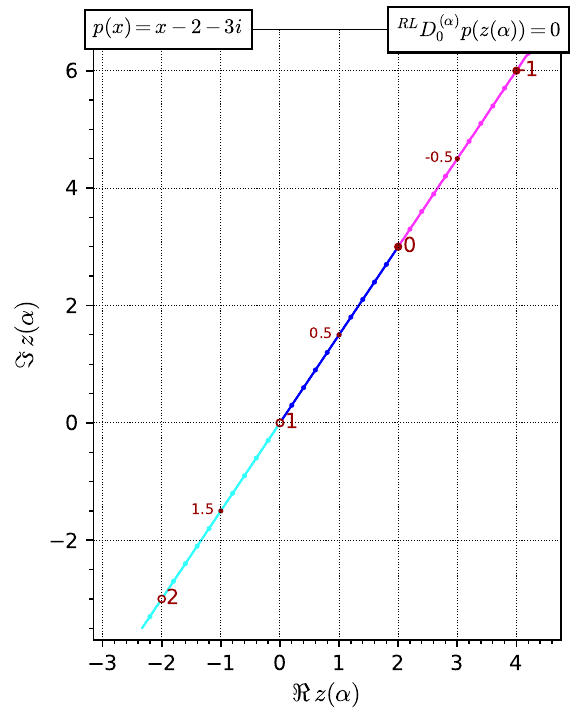}
\includegraphics[width=0.45\textwidth]{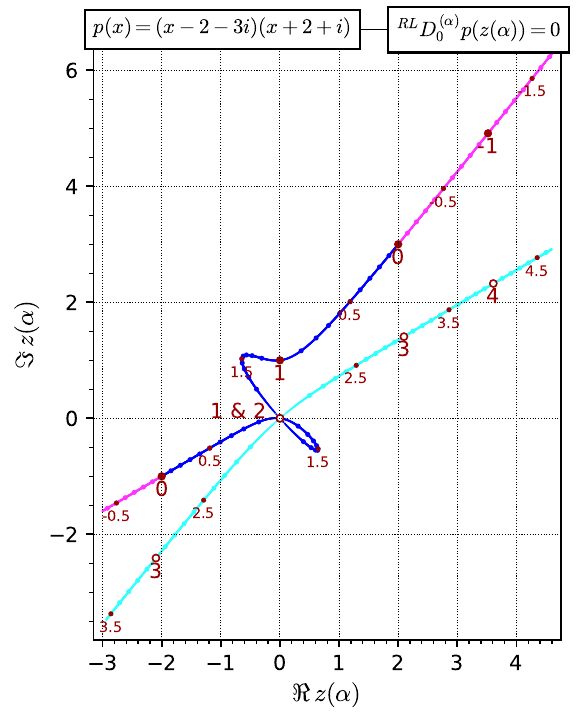}
\caption{Path $z(\alpha)$ of zeros of differintegrals of the polynomials \(p(x)=x-2-3i\) and $p(x)=(x-2-3i)(x+2+i)$.}
\label{fig one}
\end{figure}

\begin{remark} A few words should be said about our plots of the implicit functions \(z:\R\to\C\)
with \(z(\alpha)\) given by \(\RL{a}{\alpha}p(z(\alpha))=0\).  
The dots labeled `\(\bullet\)\textsf{k}' represent zeros of the \(k\)th Riemann-Liouville differintegral (thus `\(\bullet \)\textsf{0}' 
represent the zeros of the polynomial \(p(x)\) itself), while circles
`\(\circ\)\textsf{k}' represent points that are limits of \(z(\alpha)\) as \({\alpha\to k}\) but are not zeros of \(\RL{a}{k}p\).  These occur for integral \(k\) with \(k\ge\deg p(x)\), where \(\RL{a}{k}p(x)\) is constant, for example at \(x=a\), see Equation (\ref{eq at a}).
Moreover, when a point is either a zero or a limit point of zeros of both the \(j\)th and the \(k\)th differintegrals, then it is represented by `\(\bullet\) \textsf{j \& k}' or `\(\circ\) \textsf{j \& k}' respectively.  
In Figures \ref{fig two} and \ref{fig longest} we let all paths of zeros of \(\RL{a}{\alpha}p(x)\) end when they reach the origin \(a\).  In Figures \ref{fig one},  \ref{fig double}, \ref{fig p1},
\ref{fig quintic paths}, and \ref{fig quintic abs}
we continue the paths past \(a\). 
In Figures \ref{fig one}, \ref{fig double}, and \ref{fig p1} we display the path of zeros \(z(\alpha)\) for \(\alpha<0\) and \(\alpha>n\) in lighter colors than for \(0\le \alpha < n \) where \(n\) is the degree of the polynomial.  
The point \(a\) is only labeled with values for \(\alpha\ge 0\).
\end{remark}

\section{Low Degrees}\label{sec low}
Formulas for finding zeros of polynomials of low degrees have been known for centuries. Applying the Riemann-Liouville derivative to these low-degree cases
has proved to be a simple but informative exercise. In this section we summarize some of these results, stating the most useful ones as lemmas. They are 
examples of a 
dynamic that shares certain key characteristic with most high-degree cases, but some aspects of which are often unique. For example, in the linear case, 
the path the zeros take is also linear, while already in the quadratic case one observes a considerably more complex behavior.

Let us start with the linear polynomials. Here the situation is simple. The paths of zeros will always be linear, too, and they can be completely described. We have: 

\begin{lemma}[Linear Polynomials]\label{lem linear}
Let \(p(x)=(x-a)+c_0\). As $\alpha$ increases from $0$ to $1$, the path of zeros of \(\RL{a}{\alpha}p(x)\) is given by \(z(\alpha)=(\alpha-1)c_0+a\).
\begin{proof}
From (\ref{eq lem derv}), with $\beta = 1$, we get: \(\RL{a}{\alpha}p(x)=(x-a)^{-\alpha}\left((x-a)+(1-\alpha)c_0\right)\).
\end{proof}
\end{lemma}

\begin{remark}
    In addition to considering the $\alpha$-th derivatives in the usual range $0 \leq \alpha \leq n$, one could also look at what happens when $\alpha < 0$ and when $\alpha>n$.   
    In the linear case this is, again, simple.   From Lemma \ref{lem linear} we can deduce that, as with $\alpha<0$, the line of roots of the derivatives continues. Similarly, for $\alpha>1$, $\alpha\not\in\Z$ lie on the same line, see Figure \ref{fig one} below. 
\end{remark}

Let us now consider the quadratic case. This case is considerably more interesting, since there are
now two paths of zeros of the fractional derivatives, and they exhibit a much more complex and intricate behavior.

We first notice that the path of the zero closest to $a$ will directly connect with $a$, while the path of the farther zero will in the process of reaching $a$ pass through the zero of the first derivative. 
%In particular, we have: 

\begin{lemma}\label{lem quad root order}
Let $p(x)=x^2+bx+c\in\mathbb{C}[x]$ with roots \(s_1\) and \(s_2\) and let
\(d=1-\frac{4c}{b^2}\).  

For \(d\in\C\setminus\R^{<0}\), denote by \( \sqrt{d}\) the 
complex number $r$ with \(r^2=d\) and \(\Re(r)> 0\) and for $x\in\R$
let $\sign(x)=\frac{x}{|x|}$.
\begin{enumerate}
\item  If \(d\in\C\setminus\R^{<0}\) and $b\ne 0$, then for the 
roots \(s_1=\frac{-b+\sign(\Re(b))\sqrt{b^2-4c}}{2}\) and \(s_2=\frac{-b-\sign(\Re(b))\sqrt{b^2-4c}}{2}\)
of \(p(x)\) we have \(|s_2|\ge |s_1|\).
\item If $b=0$, then $|s_1|=|s_2|$.
\item If \(d\in\R^{<0}\), then \(|s_1|=|s_2|\).
\end{enumerate}
\end{lemma}

\begin{proof}
\begin{enumerate}
\item 
The roots of $p(x)$ are
\begin{align*}
\frac{-b\pm \sqrt{b^2-4c}}{2}
&=\frac{b}{2}\left(-1\pm \frac{\sqrt{b^2-4c}}{b}\right)\\
&=\frac{b}{2}\left(-1\pm \frac{\sqrt{b^2-4c}}{\sign(\Re(b))\sqrt{b^2}}\right)\\
&=\frac{b}{2}\left(-1\pm \sign(\Re(b))\sqrt{1-\frac{4c}{b^2}}\right)\\
&=\frac{b}{2}\left(-1\pm \sign(\Re(b))\sqrt{d}\right)
\end{align*}
Considering the absolute value of the last term we get:
\begin{align*}
\left|-1 \pm \sign(\Re(b))\sqrt{d}\right|^2 
&= \left(-1\pm \sign(\Re(b)) \sqrt{d}\right)\left(-1\pm \sign(\Re(b))\overline{ \sqrt{d}}\right)\\
&=1 \mp\sign(\Re(b)) \sqrt{d}\mp\sign(\Re(b))\overline{\sqrt{d}}+|d| \\
&=1 \mp2\sign(\Re(b))\Re(\sqrt{d})+|d|
\end{align*}
Because \(\Re\left( \sqrt{d}\right)> 0 \) we have
{
\begin{align*}
|-b+ \sqrt{b^2-4c}| & \ge |-b- \sqrt{b^2-4c}| &\text{ when } \sign(\Re(b))=-1 \\
|-b- \sqrt{b^2-4c}| & \ge |-b+ \sqrt{b^2-4c}| &\text{ when } \sign(\Re(b))=+1
\end{align*}
and}
\[
|-b-\sign(\Re(b)) \sqrt{b^2-4c}|  \ge |-b+\sign(\Re(b)) \sqrt{b^2-4c}|
\]
which implies $|s_2|\ge|s_1|$.
\item When $b=0$ the roots of $p(x)$ are
 $s_{1,2}=\frac{\pm \sqrt{-4c}}{2}$.  Hence
$|s_1|=|s_2|$.
\item Here $\Re\sqrt{d}=0$ and thus
$\left|-1+ \sqrt{d}\right|^2  =1+|d| =\left|-1- \sqrt{d}\right|^2$,
which implies \(|s_1|= |s_2|\).\qedhere
\end{enumerate}
\end{proof}

\begin{remark}
    As of yet, we do not know a reliable ordering of the zeros for any of the higher degree polynomials. In fact, there exist examples of cubic polynomials for which the standard Euclidean distance (which works so well for the linear and the quadratic cases) can be shown to fail: see Figure \ref{fig longest}.  
\end{remark}

 In addition to the natural ordering on the quadratic roots, another question that seems to be of interest is the one that concerns the trends of descent of their paths, especially since it had such a nice answer in the linear case. As it turns out, the asymptotes of the two quadratic paths exist, and the quadratic formula alone is enough to help us find them.

\begin{theorem}\label{prop roots quad}
For the quadratic polynomial $p(x)=(x-a)^2+c_1(x-a)+c_0$, the paths of zeros of the fractional derivatives \(\RL{a}{\alpha} p(x)\) are given as
\begin{equation}\label{eq roots quad}
z_{1,2}(\alpha)=a+\frac{-(2-\alpha) c_1\pm \sign(\Re(c_1))\sqrt{(2-\alpha)^2\cdot (c_1)^2 - 8(2-\alpha)(1-\alpha) \cdot c_0}}{4},
\end{equation}
with \(|z_2(\alpha)-a|\ge|z_1(\alpha)-a|\), for \(\alpha\in[0,2)\), and \(\lim_{\alpha\to 1}z_1(\alpha)=a\)
and 
\(\lim_{\alpha\to 2}z_{1,2}(\alpha)=a\).
\begin{proof}
With the help of (\ref{eq lem derv}), the fractional derivatives of \(p(x)\) can be written as 
\[
\RL{a}{\alpha}p(x)=\frac{2}{\Gamma(3-\alpha)}(x-a)^{-\alpha}
\left(
(x-a)^{2} +\frac{(2-\alpha)\cdot c_1}{2} (x-a) +\frac{(2-\alpha)(1-\alpha)\cdot c_0}{2}\right).
\]
Now, set $y=x-a$. Then the roots of $y^{2} +\frac{(2-\alpha)\cdot c_1}{2} y +\frac{(2-\alpha)(1-\alpha)\cdot c_0}{2}$
are 
\[
z_{1,2}(\alpha)=a+\frac{-(2-\alpha) c_1\pm \sign(\Re(c_1))\sqrt{(2-\alpha)^2\cdot (c_1)^2 - 8(2-\alpha)(1-\alpha) \cdot c_0}}{4}
\]
The ordering of the roots \(|z_1(\alpha)|\ge|z_2(\alpha)|\), for {\(\alpha\in[0,2)\)},  follows with Lemma \ref{lem quad root order}. Furthermore, we have 
\[
\lim_{\alpha\to 1} z_{1,2}(\alpha)=a+\frac{- c_1\pm \sign(\Re(c_1))\sqrt{ (c_1)^2}}{4}=a+\frac{-c_1\pm c_1}{4}.
\]
Thus $\lim_{\alpha\to 1} z_1(\alpha)=a$ and $\lim_{\alpha\to 1} z_2(\alpha)= a- \frac{c_1}{2}$.
Similarly $\lim_{\alpha\to 2}z_{1,2}(\alpha)=a$.

\end{proof}
\end{theorem}

\begin{corollary}
For \(\alpha \to \pm\infty\), the asymptotes of the quadratic paths are
\[\textstyle z_{1,2}(\alpha)\approx a+\frac{(2-\alpha) c_1}{4}\left[-1\pm \sqrt{1-\frac{8c_0}{(c_1)^2}}\right]. \]
\end{corollary}

\begin{proof}
With the quadratic formula we get
%(\ref{eq roots quad}) 
%we have,   
\begin{align*}
        z_{1,2}(\alpha)&=a+\frac{-(2-\alpha) c_1\pm(2-\alpha) c_1 \sqrt{1 - 8\frac{(1-\alpha)}{(2-\alpha)} \cdot \frac{c_0}{(c_1)^2}}}{4}\\
        &=a+\frac{(2-\alpha) c_1}{4}\left[-1\pm \sqrt{1-\frac{8c_0(1-\alpha)}{(c_1)^2(2-\alpha)}}\right],\\
\end{align*}
and since $\frac{1-\alpha}{2-\alpha}\to 1$, as $\alpha\to \pm\infty$, this yields linear asymptotes for $z_1(\alpha)$ and $z_2(\alpha)$.
\end{proof}

\begin{figure}
\includegraphics[width=0.45\textwidth]{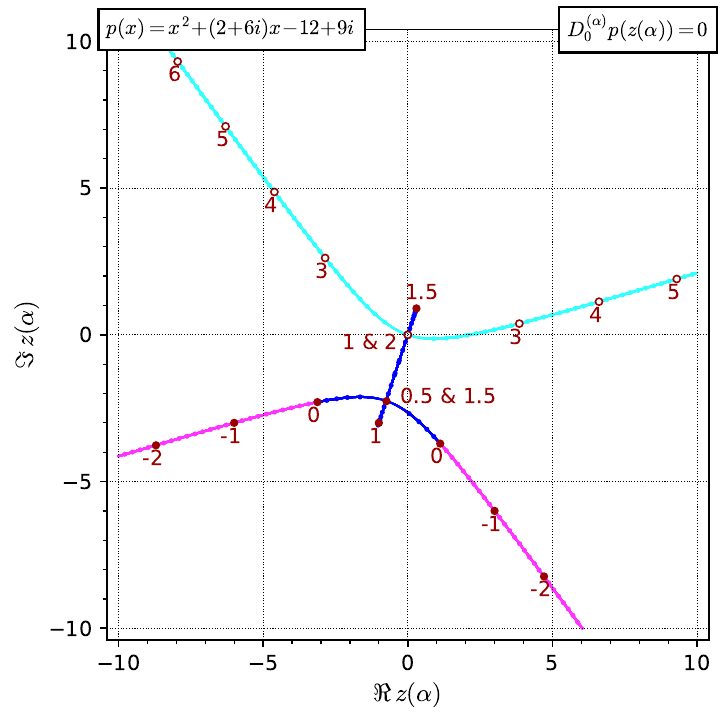}
\includegraphics[width=0.45\textwidth]{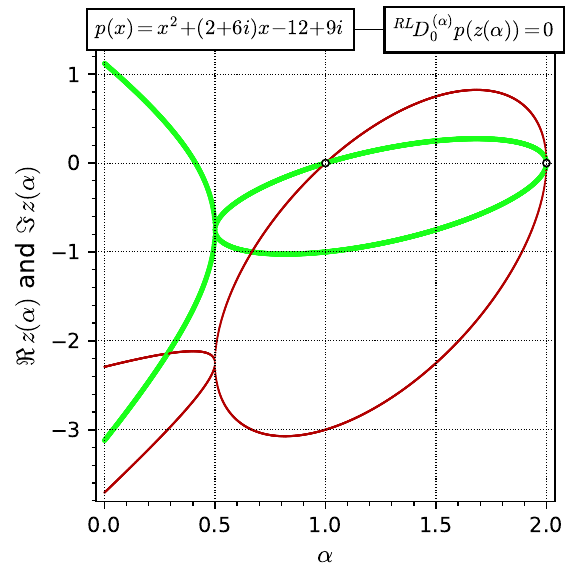}
\caption{Paths \(z(\alpha)\) of zeros of the fractional derivatives of 
\(p(x)=x^2+(2+6i)x-12+9i\) where \(\RL{0}{0.5}p\) has a double root.  In the plot on the right the wide, light green graph represents the real part of \(z(\alpha\)), while the thin, dark red graph represents its imaginary part. 
}\label{fig double}
\end{figure}

A noteworthy special case occurs when the polynomial has a double root. Then the paths display an interesting symmetry, see Figure \ref{fig double}.
In fact, it is easy to see that specializing our Theorem \ref{prop roots quad} to the case of a double zero of the polynomial itself yields:

\begin{corollary}
If $p(x)=(x-z_0)^2$ then the zeros of \(\RL{0}{\alpha}p\) are
\[z_{1,2}(\alpha)=\frac{z_0\left(-(2-\alpha) \pm  \sqrt{\alpha(2-\alpha)}\right)}{2}. \]
\end{corollary}

Another natural question to ask is whether, given that a quadratic polynomial has  distinct zeros, can its fractional derivative have a double zero. Setting \(z_1(\alpha)=z_2(\alpha)\) one gets:

\begin{corollary}\label{lem double zero}
Let $p(x)=(x-a)^2+c_1(x-a)^1+c_0(x-a)^0$. Then the fractional derivative \(\RL{a}{\alpha}(p(x))\) will have a double zero precisely for one \(\alpha\in\R\setminus\N\), namely: 
\(\alpha = 1-\frac{c_1^2}{8c_0- c_1^2}\).
%0.5=\frac{c_1^2}{8c_2- c_1^2}
%0.5(8c_2-c_1^2)=c_1^2}
%4c_2-0.5c_1^2=c_1^2
%4c_2=1.5c_1^2
%c_2=3/8 c_1^2 
\end{corollary}

\section{Flow of Zeros}\label{sec flow}

\begin{figure}
\includegraphics[width=0.45\textwidth]{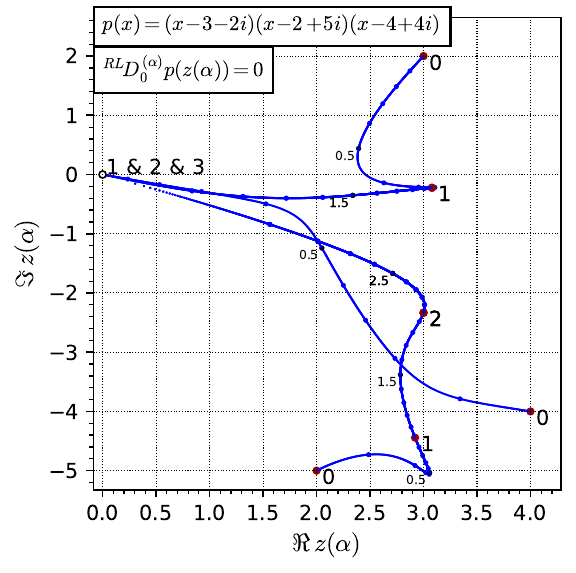}
\includegraphics[width=0.45\textwidth]{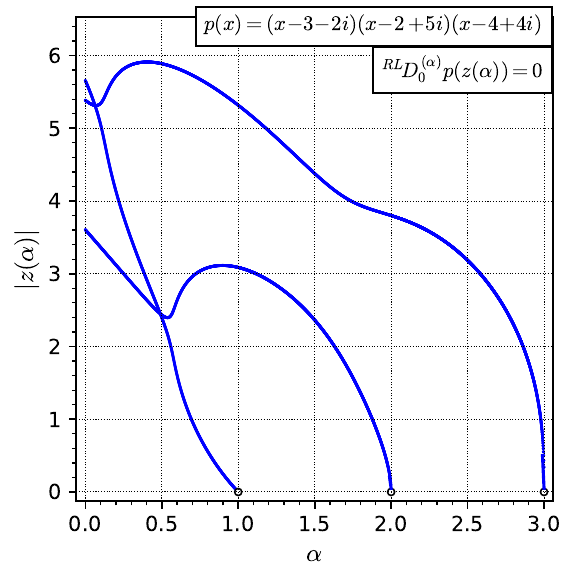}
\caption{The paths \(z(\alpha)\) of zeros of the fractional derivatives of 
the cubic \(p(x)=(x-3-2i)(x-2+5i)(x-4+4i)\)
and the absolute values \(|z(\alpha)|.\)
The zero of \(p(x)\) furthest from the origin 0 of the Riemann-Liouville Fractional derivative is
not the starting point of the longest path.
}
\label{fig longest}
\end{figure}

As stated above, one of our main goals was to consider the paths of zeros of the fractional derivatives of polynomials \(p(x) \in\C[x]\) of arbitrary degrees. Unfortunately, unlike in the linear and quadratic cases, already for the cubics we find that the situation becomes considerably more complicated. This can be seen from the fact that one of the nicest properties -- the natural ordering of zeros -- fails already for degree 3: in other words, 
it is not true in general that zeros furthest away from the origin will yield the longest paths of zeros of fractional derivatives on its way to the origin. Figure
\ref{fig longest} shows a notable counterexample.

However, certain convergence properties of the paths can be established in general. For example, the following theorem shows that all the paths terminate in the origin $a$. 

\begin{theorem}\label{theo roots lim}
Let \(p\in \mathbb{C}[x]\) of degree \(n\) such that for all \(\alpha\in [0,n)\)
the fractional derivative \(\RL{a}{\alpha}p\) has no double zeros.
Let  \(z_1,\,z_2,\,...,\,z_n\) be the zeros of $p$ and for $1\le j\le n$ let the implicit functions \(z_j:[0,n)\to \C\) given by
\(
\RL{a}{\alpha}(z_j(\alpha))=0
\)
with \(z_j(0) = z_j\).
Then there is an ordering of the roots \(z_1,\,z_2,\,...,\,z_n\) such that \(\lim_{\alpha \to j} z_j(\alpha) = a.\)
\end{theorem}

\begin{proof}
We  denote the coefficients in the expansion proved in Lemma \ref{lem derv} by
\begin{equation}
d^\alpha_j:=\prod_{k=j+1}^n\!(k-\alpha)\cdot \frac{j!}{n!}. \label{eq dalphaj}
\end{equation}
Let \(0\le j\le n\) and \(m>j\).  Here, clearly
\[
\lim_{\alpha\to m} 
d^\alpha_j=
\lim_{\alpha\to m} 
\prod_{k=j+1}^n\!(k-\alpha)\cdot \frac{j!}{n!}=0. 
\]
Denote the coefficients of the derivatives  \(\RL{a}{\alpha}p\) by
\(d^\alpha_j\,c_j\) 
as symmetric functions of the roots 
\(z_1(\alpha),z_2(\alpha), \dots, z_n(\alpha)\) of  \(\RL{a}{\alpha}p\).
Because 
\[
0=\lim_{\alpha\to m} d^\alpha_0\,c_0= \lim_{\alpha\to m} \prod_{k=1}^n z_k(\alpha)
\]
we have
\(\lim_{\alpha\to m} z_{k_1}(\alpha)=0\), for at least one \(1\le k_1\le n\).
Inductively, continuing with the next coefficient we get:
\begin{equation}\label{eq sym 1}
0=\lim_{\alpha\to m} d^\alpha_1\,c_1= \lim_{\alpha\to m} 
\sum_{l=1}^n \prod_{k\ne l} z_k(\alpha).
\end{equation}
There is one summand in (\ref{eq sym 1}) that does not contain \(z_{k_1}(\alpha)\).
So we need to have \(\lim_{\alpha\to m} z_{k_2}(\alpha)=0\) for at least one \(k_2\ne k_1\).
This argument also holds for all \(d^\alpha_j\,c_j\) with \(j<m\).
Therefore, we get \(\lim_{\alpha\to m} z_{k}(\alpha)=0\) for at least \(m\) distinct 
\(k\in\{1,\dots, n\}\).
\end{proof}

\begin{figure}
\includegraphics[width=0.45\textwidth]{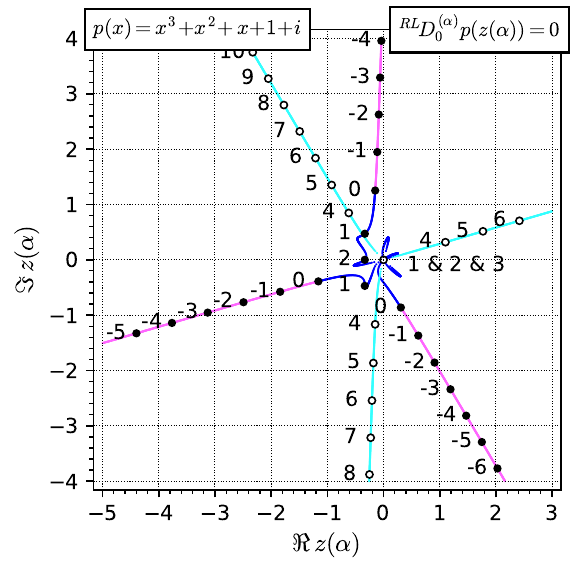}
\includegraphics[width=0.45\textwidth]{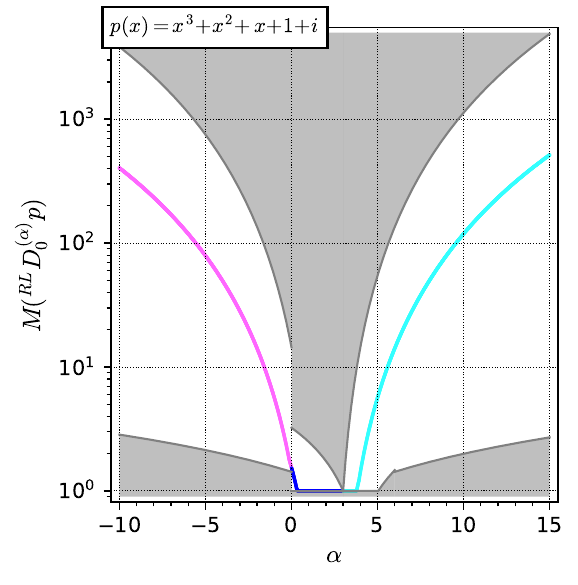}
\caption{Paths \(z(\alpha)\) of zeros of \(\RL{0}{\alpha}(x^3+x^2+x+1+i)\).
illustrating the growth of the absolute value of \(z(\alpha)\) as
\(\alpha\to\infty\) and \(\alpha\to-\infty\) and the Mahler measure
of \(\RL{0}{\alpha}p\) along with the bounds from Theorems \ref{theo upper bound} and \ref{theo lower bound}.
The missing points for $\alpha\in\N$, $n\ge 4$ are not shown on the right.}
\label{fig p1}
\end{figure}

\section{Bounds}\label{sec bounds}

As stated in the introduction, for \(f\in\C[x]\) the Gauss-Lucas theorem states that all zeros of \(f'\) lie in the convex hull of the set of zeros of \(f\), see \cite[Theorem 6.1]{marden1966}. By induction this generalizes to all integral derivatives. Unfortunately, although all roots of the fractional derivatives converge to the origin, by our Theorem \ref{theo roots lim}, the analogue of the Gauss-Lucas theorem does not hold for the fractional derivatives.
This is an immediate consequence of a result by Genchev and Sendov
{\cite{genchev-sendov}}, which is also stated as {\cite[Theorem 2]{nikolov-sendov}}:

\begin{theorem}
Let \(L:\C[x]\to\C[x]\) be a linear operator, such that \(L(p)\ne 0\) implies that
the convex hull of the set of roots of \(p\) contains the roots of \(L(p)\). 
Then \(L\) is a linear functional or there are \(c\in\C\setminus\{0\}\) and \(k\in\N\) such that \(L(p)=cp^{(k)}\).
\end{theorem}

Figure \ref{fig longest} illustrates this result by giving a specific counterexample to the Gauss-Lucas property for the case of the Riemann-Liouville fractional derivatives. Nevertheless, it is possible to make some useful statements about how the absolute values of zeros
\(z_k(\alpha)\)
of the fractional derivatives \(\RL{0}{\alpha} f\)
decrease as \(\alpha\) increases in terms of the Mahler measure of \(f\).

Let \(p(x)=x^n+\sum_{j=1}^{n-1} c_j x^j=\prod_{j=1}^n (x-z_j)\in\mathbb{C}[x]\).
For the Mahler measure \(M(f)\) \cite{mahler} we have
\[
M(p) = \exp\left(\int_0^1 \log(|f(e^{2\pi i\theta})|)\, d\theta \right)= \prod_{j=1}^n \max\{1,|z_j|\}.
\]
Denote the height of \(f\) by
\(
\nheight{p} = \max\{c_0,\dots,c_n\}
\)
and the length of \(f\) by
\(
\nlength{p}=|c_0|+\dots+|c_n|.
\)
Recall that Mahler was able to prove the bounds 
\begin{equation}\label{eq mahler height}
\binom{n}{\lfloor n/2 \rfloor}^{-1} \nheight{p} \le M(p) \le \nheight{p}  \sqrt{n+1}
\end{equation}
and
\begin{equation}\label{eq mahler length}
2^{-n} \nlength{p} \le M(p) \le  \nlength{p}.
\end{equation}
For \(\ntwo{p} = \left(\sum_{j=1}^n |c_j|^2\right)^\frac{1}{2}\) we have Landau's inequality
\cite{landau1905}
\begin{equation}\label{eq landau}
M(p)\le ||p||_2.
\end{equation}

We generalize the definition of the Mahler measure to fractional derivatives.
Let \(Z(\alpha)\) be the set of zeros of \(\RL{0}{\alpha}p\).  We set
\[
M\left(\RL{0}{\alpha}p\right)=\prod_{z\in Z(\alpha)}\max\{1,|z|\}.
\]
and prove bounds similar to (\ref{eq mahler height}),
(\ref{eq mahler length}), and (\ref{eq landau})
for the fractional cases. 
We first estimate the coefficients \(d^\alpha_j\) from Proposition \ref{theo roots lim} and then use them to derive bounds for \(M\left(\RL{0}{\alpha}p\right)\).

\begin{figure}
\includegraphics[width=0.45\textwidth]{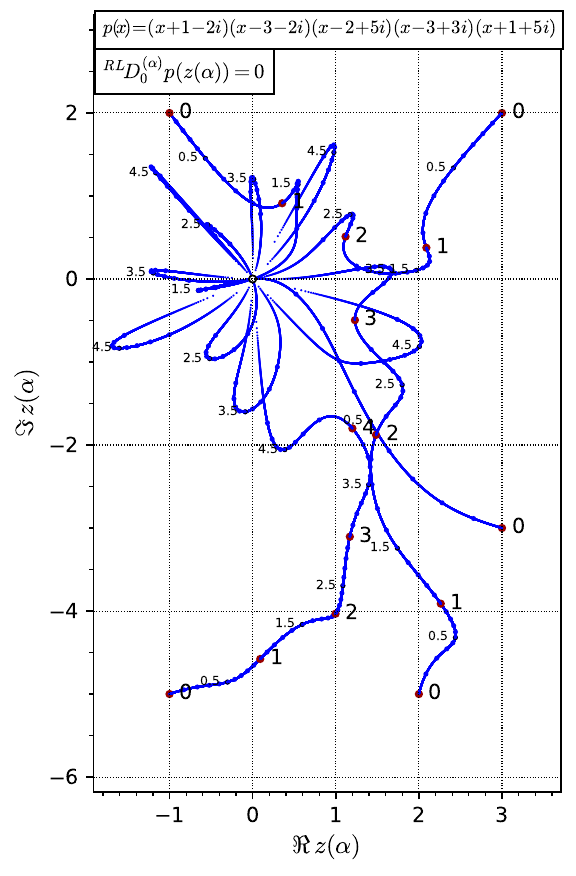}
\includegraphics[width=0.45\textwidth]{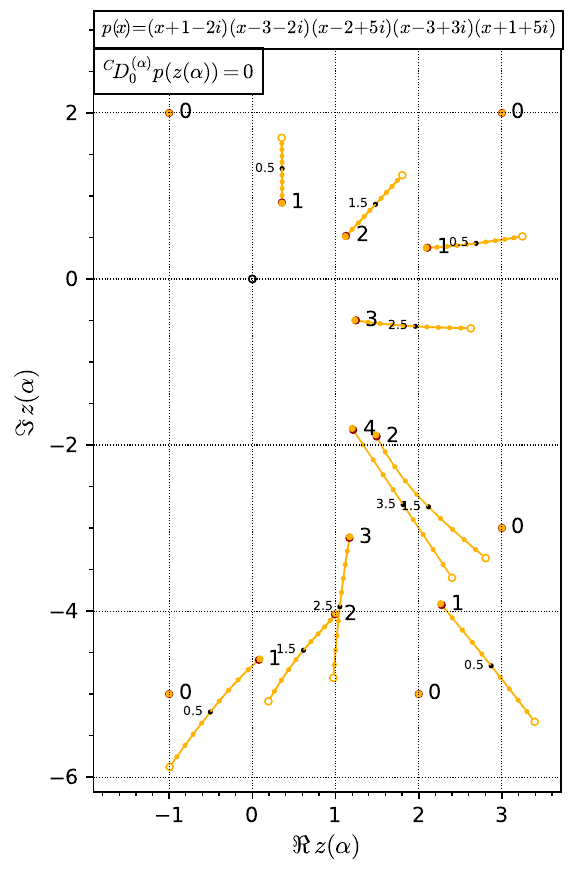}
\caption{Comparing paths \(z(\alpha)\) of zeros of the Riemann-Liouville (left) and Caputo (right) fractional derivatives 
of the  quintic 
\(p(x)=(x+1-2i)(x-3-2i)(x-2+5i)(x-3+3i)(x+1+5i)\), for \(0\le\alpha\le 5\).
}
\label{fig quintic paths}
\end{figure}

\begin{lemma}\label{lem dj}
Let \(n\in\N\) and \(\alpha\in\R\setminus\N\).
Let
\(d^\alpha_j=\left(\prod_{k=j+1}^n\!(k-\alpha)\right)\cdot \frac{j!}{n!}
\) where \(0 \le j\le n-1\).  Then
\begin{enumerate}
\item\label{dj gt 0}
\(|d^\alpha_j| \le \frac{n-\alpha}{n}\), for  \(0<\alpha< n\).
\item\label{dj lt}
\(|d^\alpha_j| \le \frac{\ff{|n-\alpha|}{n}}{n!}\), 
where \(\ff{|n-\alpha|}{n}=\prod_{k=0}^{n-1} |n-\alpha-k|\), % denotes the falling factorial
for \(\alpha<0\) and \(\alpha>n\).
\item\label{dj lt 0 gt 2n}
\( |d^\alpha_j|\ge\frac{|n-\alpha|}{n}\),  for \(\alpha<0\)
and \(\alpha>2n\).
\item\label{dj gt n}
\(|d^\alpha_j|\ge\frac{\alpha-n}{n} \binom{n-1}{\lceil\frac{n-1}{2}\rceil}^{-1}\), for \(n<\alpha\le 2n\).
\end{enumerate}
\end{lemma}
\begin{proof}
{
\begin{enumerate}
\item
For \(0\le\alpha < j+1\), we have
\[
|d^\alpha_j| = \frac{n-\alpha}{n} \prod_{k=j+1}^{n-1}\frac{k-\alpha}{k} \cdot \frac{j!}{j!}< \frac{n-\alpha}{n}.
\]
When \(j+1 < \alpha\) set \(h:=\lfloor\alpha\rfloor\).  We get
\begin{align*}
|d^\alpha_j| 
&=
\frac{|n-\alpha|\cdots|h+1-\alpha|\cdot|h-\alpha|\cdots|j+1-\alpha|\cdot j!}{n!}\\
%&=
%\frac{n-\alpha}{n}\cdot\frac{(n-1-\alpha)\cdots(h+1-\alpha)\cdot(\alpha-j-1)\cdots(\alpha-h)\cdot j!}{(n-1)!}\\
%&\le
%\frac{n-\alpha}{n}\cdot\frac{(n-1-h)\cdots 1 \cdot(h-j)\cdots 1\cdot j!}{(n-1)!}\\
&\le \frac{n-\alpha}{n}\cdot\frac{(n-h-1)! \cdot (h-j)!\cdot j!}{(n-1)!}
\le \frac{n-\alpha}{n}.
\end{align*}
\item For \(\alpha<0\) and \(\alpha>n\), we have
\begin{align*}
|d^\alpha_j| 
&=
\frac{|n-\alpha|\cdots|j+1-\alpha|\cdot j!}{n!}
=
\frac{|\alpha-n|\cdots|\alpha-j-1|\cdot j!}{n!}\\
&\le
\frac{|\alpha-n|\cdots|\alpha-1|}{n!}=\frac{|n-\alpha|_n}{n!}.
\end{align*}
\item
{
For $\alpha<0$ we have
\[
|d_j^{(\alpha)}|=
\frac{n-\alpha}{n}\cdot\frac{(n-1-\alpha)\dots(j+1-\alpha)\cdot j!}{(n-1)!}\ge
%\frac{n-\alpha}{n}\cdot\frac{(n-1)\dots(j+1)\cdot j!}{(n-1)!}=
\frac{\alpha-n}{n}=\frac{|n-\alpha|}{n}.
\]
For $\alpha>2n$ we have
\[
|d_j^{(\alpha)}|=
\frac{n-\alpha}{n}\cdot\frac{(\alpha-n-1)\dots(\alpha-j+1)\cdot j!}{(n-1)!}
\ge\frac{n-\alpha}{n}=\frac{|n-\alpha|}{n}.
\]
}
\item
For \(2n>\alpha>n\), we have
\[
|d^\alpha_j| 
%&= \frac{|n-\alpha|\cdots|j+1-\alpha|\cdot j!}{n!}=\frac{(\alpha-n)\cdots(\alpha-j-1)\cdot j!}{n!}\\
%&=\frac{\alpha-n}{n}\cdot\frac{(\alpha-n+1)\cdots(\alpha-j-1)\cdot j!}{(n-1)!} 
\ge\frac{\alpha-n}{n}\cdot\frac{1\cdots(n-j-1)\cdot j!}{(n-1)!}
\ge\frac{\alpha-n}{n} \binom{n-1}{\lceil\frac{n-1}{2}\rceil}^{-1}.\quad\qedhere
\]
\end{enumerate}
}
\end{proof}

\begin{figure}
\includegraphics[width=0.45\textwidth]{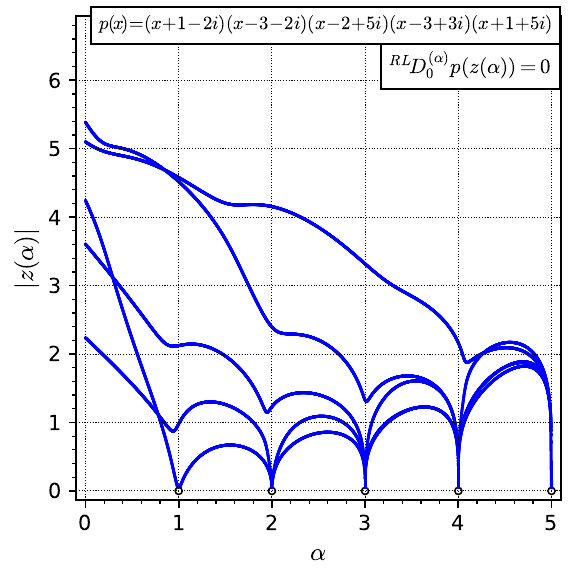}
\includegraphics[width=0.45\textwidth]{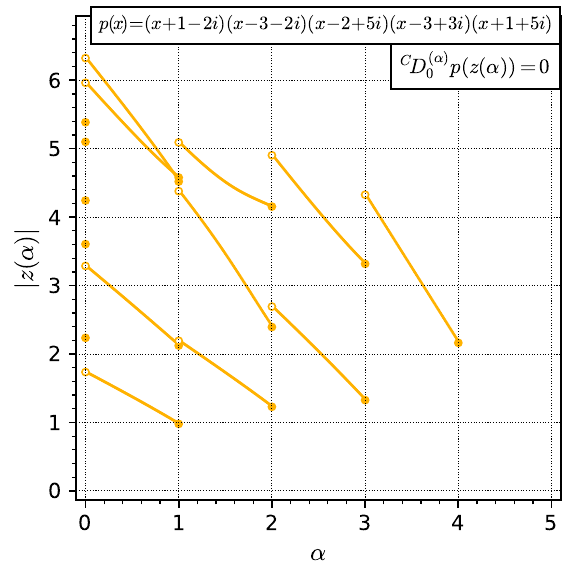}
\caption{Comparing the absolute value of \(z(\alpha)\) of zeros of the Riemann-Liouville (left) and Caputo (right) fractional derivatives 
of the  quintic 
\(p(x)=(x+1-2i)(x-3-2i)(x-2+5i)(x-3+3i)(x+1+5i)\) from Figure \ref{fig quintic paths} for \(0\le\alpha\le 5\).
}
\label{fig quintic abs}
\end{figure}

We now show that \(M\left(\RL{0}{\alpha} p\right)\) 
is bounded from above for \(0<\alpha<\deg p=n\) where the bound linearly 
decreases with \(\alpha\).  Furthermore, in Theorem \ref{theo lower bound},
we see that for \(\alpha<0\) and \(\alpha>\deg p\) the Mahler measure
\(M\left(\RL{0}{\alpha} p\right)\) increases at least linearly with \(\alpha\).

{
We use the notation from Lemma \ref{lem derv}.
For $\alpha\in\R\setminus\N^{\ge n}$ we have
\(
%\begin{align*}
M\left(\RLs{0}{\alpha}f\right)  
%=M\left(\frac{\Gamma(n+1-\alpha)}{n!}x^\alpha\cdot \RL{a}{\alpha} f\right)
= M\left(\RL{0}{\alpha}f\right)
%\end{align*}
\)
since \(\RLs{0}{\alpha}f\) and \(\RL{0}{\alpha}f\) have the same roots except for an additional root \(0\) in some cases.
}
\begin{theorem}\label{theo upper bound}
Let \(f(x)\in\mathbb{C}[x]\) be monic and let \(0<\alpha<n\).  Then
\begin{enumerate}
\item $M(\RL{0}{\alpha} f)  \le \frac{n-\alpha}{n}\nlength{f}+1$
\item $M(\RL{0}{\alpha} f)  \le  \sqrt{n+1} \max\left\{\frac{n-\alpha}{n} \nheight{f},1\right\}$
\item $M(\RL{0}{\alpha} f)  \le \frac{n-\alpha}{n} \ntwo{f}+1$
\end{enumerate}
\end{theorem}
%\fixme{In some places (first line of proof) we are implicitly assuming that \(f\) is monic.  Shall we explicitly say \(f\) is monic or remove the implicit assumptions ?}
\begin{proof}
\begin{enumerate}
\item
{
With Lemma \ref{lem dj} (\ref{dj gt 0}), we get
\begin{align*}
\nlength{\RLs{0}{\alpha}}
&=|d_0^{(\alpha)}c_0|+|d_1^{(\alpha)}c_1|+\dots+|d_{n-1}^{(\alpha)}c_{n-1}|+1\\
&\le \frac{n-\alpha}{n}\left(|c_0|+|c_1|+\dots+|c_{n-1}|\right)+1\\
&\le \frac{n-\alpha}{n} \nlength{f}+1.
\end{align*}
}
With Mahler's (\ref{eq mahler length}), we get
\[
M\left(\RL{a}{\alpha} f\right) =
M\left(\RLs{0}{\alpha}\right) \le
\nlength{\RLs{0}{\alpha}} \le \frac{n-\alpha}{n}\nlength{f}+1.
\]
\item
{
By Lemma \ref{lem dj} (\ref{dj gt 0}), we have
\begin{align*}
\nheight{\RLs{0}{\alpha}f}
&=\max\{|d_0^{(\alpha)}c_0|,|d_1^{(\alpha)}c_1|,\dots,|d_{n-1}^{(\alpha)}c_{n-1}|,1\}\\
&\le \max\left\{\frac{n-\alpha}{n}\max\{|c_0|,|c_1|,\dots,|c_{n-1}|\},1\right\}
\\
%&\le \max\left\{\frac{n-\alpha}{n}\max\{|c_0|,|c_1|,\dots,|c_{n-1}|,1\},1\right\}
%\\
&\le\max\left\{\frac{n-\alpha}{n}\nheight{f},1\right\}\\
%&\le \frac{n-\alpha}{n}\nheight{f}+1.
\end{align*}
}
With Mahler's (\ref{eq mahler height}),
we get
\[
M\left(\RL{0}{\alpha} f\right) =
M\left(\RLs{0}{\alpha} f\right) \le
 \sqrt{n+1} \nheight{\RLs{0}{\alpha}} \le  \sqrt{n+1}  \max\left\{\frac{n-\alpha}{n} \nheight{f},1\right\}.
\]
\item
By Lemma \ref{lem dj} (\ref{dj gt 0}), we have
\begin{align*}
\ntwo{\RLs{0}{\alpha}}f
&= \sqrt{(d_0^{(\alpha)}c_0)^2+\dots+(d_{n-1}^{(\alpha)}c_{n-1})^2+1}\\
%&\le \sqrt{\left(\frac{n-\alpha}{n}c_0\right)^2+\dots+
%\left(\frac{n-\alpha}{n}c_{n-1}\right)^2+1}\\
%&\le \sqrt{\left(\frac{n-\alpha}{n}c_0\right)^2+\dots+
%\left(\frac{n-\alpha}{n}c_{n-1}\right)^2}+1\\
&\le\frac{n-\alpha}{n} \sqrt{c_0^2+\dots+c_{n-1}^2}+1\\
%&\le\frac{n-\alpha}{n} \sqrt{c_0^2+\dots+c_{n-1}^2+1}+1\\
&\le\frac{n-\alpha}{n}\ntwo{f}+1.
\end{align*}
With Landau's (\ref{eq landau}), we get
\[
M\left(\RL{0}{\alpha} f \right)=
M\left(\RLs{0}{\alpha}\right) \le \ntwo{\RLs{0}{\alpha}} \le\frac{n-\alpha}{n} ||f||_2+1.
\qedhere\]
\end{enumerate}
\end{proof}

For the remaining values of $\alpha$ we get:
%Now we are able to prove:

\begin{theorem}\label{theo lower bound}
Let \(f(x)\in\C[x]\) be monic of degree \(n\) and let $\alpha\in\R\setminus\N$.  Then
\begin{enumerate}
\item\label{tlb1}\(M(\RL{0}{\alpha} f)\ge 2^{-n}\left(1+\frac{|n-\alpha|}{n}(\nlength{f}-1)\right)\), for \(\alpha<0\) and \(\alpha>2n\)
\item\label{tlb2}\(M(\RL{0}{\alpha} f)\le \frac{\ff{|n-\alpha|}{n}}{n!}\nlength{f}+1\), for  \(\alpha<0\) and \(\alpha>n\)
\item\label{tlb3}\(M(\RL{0}{\alpha} f)\ge 2^{-n}\binom{n-1}{\lceil(n-1)/2\rceil}^{-1}\frac{\alpha-n}{n}\nlength{f}\), for \(n<\alpha< 2n\)
%\item\label{tlb4}\(M(\RL{0}{\alpha} f)\le 1\) for \(n<\alpha\le 2n\)
\end{enumerate}
\begin{proof}
\begin{enumerate}
\item[\ref{tlb1}.]
By Lemma \ref{lem dj} (\ref{dj lt 0 gt 2n}), we have 
\begin{align*}
\nlength{\RLs{0}{\alpha}f}
&=|d_0^{(\alpha)}c_0|+|d_1^{(\alpha)}c_1|+\dots+|d_{n-1}^{(\alpha)}c_{n-1}|+1\\
&\ge \frac{|n-\alpha|}{n}\left(|c_0|+|c_1|+\dots+|c_{n-1}|\right)+1\\
%\left(|c_0|+|c_1|+\dots+|c_{n-1}|+1-1\right)+1\\
%&\ge \frac{|n-\alpha|}{n}\left(|c_0|+|c_1|+\dots+|c_{n-1}|+1-1\right)+1\\
&=1+\frac{|n-\alpha|}{n}(\nlength{f}-1).
\end{align*}
With Mahler's (\ref{eq mahler length}), we get
\[
M(\RL{0}{\alpha} f) =M\left(\RLs{0}{\alpha}\right) \ge 2^{-n} \nlength{\RLs{0}{\alpha}} \ge 2^{-n}\left(1+\frac{|n-\alpha|}{n}(\nlength{f}-1)\right).
\]
\item[\ref{tlb2}.]
{
By Lemma \ref{lem dj} (\ref{dj lt}), we have
\begin{align*}
\nlength{\RLs{0}{\alpha}f}
&= \sqrt{(d_0^{(\alpha)}c_0)^2+\dots+(d_{n-1}^{(\alpha)}c_{n-1})^2+1}\\
&\le\frac{\ff{|n-\alpha|}{n}}{n!}\sqrt{c_0^2+\dots+c_{n-1}^2}+1\\
&\le\frac{\ff{|n-\alpha|}{n}}{n!}\ntwo{f}+1.
%&=|d_1^{(\alpha)}c_1|+\dots+|d_{n-1}^{(\alpha)}c_{n-1}|+1\\
%&\le \frac{\ff{|n-\alpha|}{n}}{n!}\left(|c_1|+\dots+|c_{n-1}|\right)+1\\
%&\le \frac{\ff{|n-\alpha|}{n}}{n!}\nlength{f_1}+1
\end{align*}
}
%&\le \sqrt{\left(\frac{n-\alpha}{n}c_0\right)^2+\dots+
%\left(\frac{n-\alpha}{n}c_{n-1}\right)^2+1}\\
%&\le \sqrt{\left(\frac{n-\alpha}{n}c_0\right)^2+\dots+
%\left(\frac{n-\alpha}{n}c_{n-1}\right)^2}+1\\
%&\le\frac{n-\alpha}{n} \sqrt{c_0^2+\dots+c_{n-1}^2+1}+1\\
With Landau's (\ref{eq landau}),
we get
\[
M\left(\RL{a}{\alpha} f\right) =
M\left(\RLs{a}{\alpha} f\right) \le
\nlength{\RLs{a}{\alpha} f}
\le \frac{\ff{|n-\alpha|}{n}}{n!}\nlength{f}+1.
\]
\item[\ref{tlb3}.]
%Follows from  (\ref{eq mahler length}) and Lemma \ref{lem dj} (\ref{dj gt n}).
With Lemma \ref{lem dj} (\ref{dj gt n}) we get
{
\begin{align*}
\nlength{\RLs{0}{\alpha}}
&=|d_0^{(\alpha)}c_0|+|d_1^{(\alpha)}c_1|+\dots+|d_{n-1}^{(\alpha)}c_{n-1}|+1\\
&\ge \frac{\alpha-n}{n}\binom{n-1}{\lceil(n-1)/2\rceil}^{-1}\left(|c_0|+|c_1|+\dots+|c_{n-1}|+1\right)\\
&=\frac{\alpha-n}{n}\binom{n-1}{\lceil(n-1)/2\rceil}^{-1}\nlength{f}
\end{align*}
}
With Mahler's (\ref{eq mahler length}), we get
\[
M(\RL{0}{\alpha} f) =M\left(\RLs{0}{\alpha}\right) \ge 2^{-n} \nlength{\RLs{0}{\alpha}} \ge 2^{-n}\frac{\alpha-n}{n}\binom{n-1}{(n-1)/2}^{-1}\nlength{f}.\qedhere
\]
%\item[\ref{tlb4}.]
%\begin{align*}
%M(\RL{0}{\alpha}
%\nlength{\RLs{0}{\alpha}}
%&=
%|d_0^{(\alpha)}c_0|+|d_1^{(\alpha)}c_1|+\dots+|d_{n-1}^{(\alpha)}c_{n-1}|+1\\
%&=
%\end{align*}
\end{enumerate} 
\end{proof}
\end{theorem}

In Figure \ref{fig p1} we present the paths of zeros and the
Mahler measures of the fractional derivatives of a degree 3 polynomial along with the bounds from Theorems \ref{theo lower bound} and \ref{theo upper bound}.
Furthermore Figures \ref{fig one}, \ref{fig double}, and \ref{fig p1} show 
the growth of $M(\RL{0}{\alpha} p)$ for \(\alpha<0\) and \(\alpha>0\).

\begin{figure}
\begin{center}
\includegraphics[width=0.45\textwidth]{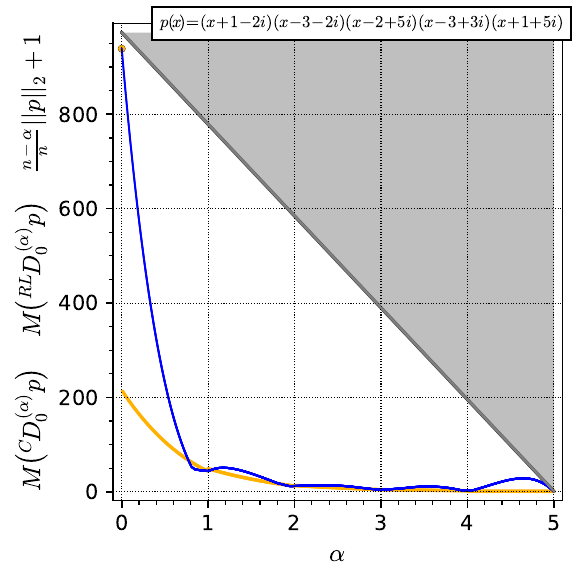}
\end{center}
\caption{Mahler measure of the Riemann-Liouville and Caputo fractional derivatives of the quintic \(p(x)=(x+1-2i)(x-3-2i)(x-2+5i)(x-3+3i)(x+1+5i)\)  for \(0\le\alpha\le 5\) along with the bound from
Theorem \ref{theo upper bound} and Corollary \ref{cor upper bound caputo}.}
\label{fig bound mahler}
\end{figure}

\section{Caputo Fractional Derivatives} \label{caputo}

Switching the order of differentiation and integration in Definition \ref{def rl} yields the Caputo fractional derivative, see \cite{diethelm2010} for example.

\begin{definition}[Caputo fractional derivative]\label{def caputo} %\cite{ortigueira}% ]
Let \(f\) be analytic on a convex open set \(C\) let and \(a\in C\).  Let $\alpha>0$ and $m=\lceil \alpha \rceil$ then the \(\alpha\)-th Caputo fractional derivative is
\[
\CA{a}{\alpha} f(t)
=
\RLI{a}{m-\alpha} 
\frac{d^m}{dt^m} 
f(t).
\]
%
%\frac{1}{\Gamma(m-\alpha)}\int_{t_0}^t (t-\tau)^{m-%\alpha-1}f^{(m)}(\tau)\,d\tau
%\]
\end{definition}
The Caputo fractional derivative has the advantage that the composition rule holds, while the power rule is very similar to that of the Riemann-Liouville fractional derivative: 
\[
    \CA{a}{\alpha}(x-a)^\beta= \begin{cases}
        0 & 
        \begin{array}{l}\text{if } \beta \in \{0,1,\dots, m-1\}
        \end{array}\\
        \dfrac{\Gamma(\beta+1)}{\Gamma(\beta-\alpha+1)}(x-a)^{\beta-\alpha} & 
\begin{array}{l}        
        \text{if } \beta\in\N \text{ and }\beta \ge  m\text{ or }\\[-0.5ex]
        \text{if } \beta\not\in\N \text{ and }\beta > m-1,
\end{array}
    \end{cases}
\]
where $m=\lceil \alpha \rceil$.

These derivatives differ, but obey some common general trends. Figures \ref{fig quintic paths} and \ref{fig quintic abs} compare the paths of zeros of the Riemann-Liouville and Caputo derivatives of a degree 5 polynomial.

The upper bounds from Theorem \ref{theo upper bound} easily transfer 
to the Caputo fractional derivative.  Because the coefficients of the Caputo fractional derivatives of a polynomial \(p\in\C[x]\) are either the 
same as those of the Riemann-Liouville fractional derivative or zero, 
we have
\[
\norm{\CAs{0}{\alpha}}{k} \le \norm{\RLs{0}{\alpha}}{k},
\]
for  \(k\in\{1,2,\infty\}\).  This yields the bounds:
\begin{corollary}
\label{cor upper bound caputo}
Let \(p\in\C[x]\) be monic of degree $n$ and let \(0<\alpha<n\).  Then 
\begin{enumerate}
\item $M(\CA{0}{\alpha} p)  \le \frac{n-\alpha}{n}\nlength{p}+1$,
\item $M(\CA{0}{\alpha} p)  \le  \sqrt{n+1} \max\left\{\frac{n-\alpha}{n} \nheight{p},1\right\}$,
\item $M(\CA{0}{\alpha} p)  \le \frac{n-\alpha}{n} \ntwo{p}+1$.
\end{enumerate}
\end{corollary}

\section{Conclusion} \label{conclusion}

The bounds we have established in Section \ref{sec bounds} and Section \ref{caputo} were sufficient for our purposes, but they are 
far from best possible. 
Figure \ref{fig bound mahler} illustrates the decline of $M(\RL{0}{\alpha} p)$ and $M(\CA{0}{\alpha} f)$, when $\alpha$ approaches \(n\), as described by Theorem \ref{theo upper bound} and Corollary \ref{cor upper bound caputo} and Theorem \ref{theo roots lim}. In a future work, it would be interesting to consider the true growth of the paths of zeros, for \(\alpha<0\) and \(\alpha>n\). Also, the maximal extent of loops of paths, after traversing the origin \(a\), is something that could be worth looking at. Moreover, as Figure \ref{fig p1} clearly shows, the paths exhibit very distinct linear asymptotes in both directions \(\alpha\to\infty\) and \(\alpha\to-\infty\). For polynomials of higher degrees, their exact directions are not yet known. 

In addition to this, as we have noted earlier, the useful Gauss-Lucas property does not hold universally for the fractional derivatives of polynomials. However, some of the dynamical properties we have observed could be investigated with insights related to those that play a key role in the integral case. In particular, Gauss himself suggested a very intriguing physical interpretation of the nontrivial critical points of a polynomial (the critical points which are not zeros) as the equilibrium points in certain force fields, generated by particles placed at the zeros of the polynomial, with masses equal to the multiplicity of the zeros and repelling with a force inversely proportional to the distance. This amazing physical application of a purely theoretical polynomial concept is exceedingly intriguing and should be investigated further. It could go a long way in explaining the profound intricacies of the paths of zeros, and their seemingly chaotic local behavior.  

\section{Acknowledgments}
We would like to thank the anonymous reviewer for the useful comments and for pointing out a couple of mistakes.
All plots were created with the computer algebra system SageMath \cite{sage}.

\bibliography{zeta}
\bibliographystyle{amsalpha}
\end{document}